\documentclass[12pt]{article}
\usepackage{bbm}
\usepackage{amsfonts}
\usepackage{lmodern}
\usepackage{pifont}
\usepackage{dashrule}
\usepackage{hyperref}
\usepackage[all,pdf]{xy}
\usepackage{mathrsfs}
\usepackage{indentfirst}
\usepackage{amsmath}
\usepackage{amssymb}
\usepackage[amsmath, thmmarks]{ntheorem}
\usepackage{cite}
\usepackage{graphicx}
\usepackage{tikz}
\newtheorem{definition}{\bf Definition}[section]

\newtheorem{theorem}{\bf Theorem}[section]
\newtheorem{remark}{\bf Remark}[section]
\newtheorem{corollary}{\bf Corollary}[section]
\newtheorem{proposition}{\bf Proposition}[section]
\newtheorem{example}{\bf Example}[section]

\begin{document}
\setcounter{page}{1}

\title{{\textbf{Nullnorms on bounded trellises}}\thanks {Supported by National Natural Science Foundation of China (No. 11871097, 12271036).}}
\author{  Zhenyu Xiu $^1$\footnote{\emph{E-mail address}: xyz198202@163.com},  Xu Zheng$^2$\footnote{Corresponding author. \emph{E-mail address}: 3026217474@qq.com  }\\
\emph{\small $^{1,2}$  College of Applied Mathematics, Chengdu University of Information Technology, }\\
\emph{\small Chengdu 610000,  China }}
\newcommand{\pp}[2]{\frac{\partial #1}{\partial #2}}
\date{}
\maketitle

\begin{quote}
{\bf Abstract}
In this paper, we introduce the notion of nullnorms on bounded trellises and study some basic properties. Based on the existence of $t$-norms and $t$-conorms on arbitrary
bounded trellises, we propose some construction methods of nullnorms on bounded trellises. Moreover, some illustrative examples are provided.

\end{quote}

\section{Introduction}\label{intro}

The nullnorms and $t$-operators  on the unit interval $[0, 1]$, as generalizations of triangular norms ($t$-norms, for short) and triangular conorms ($t$-conorms, for short), were introduced by Calvo et al. \cite{TC01} and Mas et al. \cite{MM99}, respectively. In particular, Mas et al. \cite{MM02} showed that nullnorms and $t$-operators are equivalent since they have the same block structures in $[0, 1]^{2}$. Since then, many results with respect to nullnorms are proposed by researchers \cite{JD08,PD04,PD15,FQ05} and nullnorms  have been proven to have applications in several fields, such as fuzzy logic, fuzzy sets theory, expert systems, neural networks
and so on \cite{CA06,DD85,DD00,MG09,RM98}.
Based on  the existence of $t$-norms and $t$-conorms on an arbitrary bounded lattice $L$, Kara\c{c}al et al. \cite{FK15} generalized the concept of nullnorms from the unit
interval $[0,1]$ to bounded lattices and proposed some construction methods for nullnorms on bounded lattices.
Following this article, a number of construction methods for nullnorms on bounded lattices were proposed by many researchers \cite{GC18,UE18,GC19,GC20,GC201,XH22,XS20,XW20}.

Recently, the concepts of $t$-norms had been generalized from the bounded lattices to bounded trellises by  Zedam and  De Baets \cite{LZ23} and some construction methods had been proposed.
Subsequently,  uninorms on  bounded trellises were studied by  Kong and Zhao \cite{YK24}.
The order on a  trellis is a binary relation that satisfies reflexivity, antisymmetry, but does not need to satisfy transitivity. It is obvious that  the bounded trellis is more general than the bounded lattice because the latter  needs to satisfy reflexivity, antisymmetry, transitivity.
In \cite{LZ23},   Zedam and  De Baets wrote:  many theoretical and practical developments warrant us to look beyond transitivity.
For example, in \cite{BD06} and \cite{BD15}, the study of (absence of) transitivity in
the comparison of random variables resulted in the framework of cycle-transitivity.  In fact,
the  absence of transitivity can lead to  the presence of cycles or preference loops
 or simply incomparability. For instance, cycles appear in species competition structures such
as tournaments preventing extinction and supporting biodiversity \cite{BK02,TR07}, usually catalogued under the Rock-Paper-Scissors metaphor and  incomparability always exists, such as the intransitivity of indifference \cite{PF70}.
 Moreover, in  \cite{YK24}, more nontransitive relations were founded  in games, the concept of closeness,  graph theory and logic of nontransitive implications.

As  we know,  the concepts of  $t$-norms ($t$-conorms) \cite{SS06,GD181}, uninorms \cite{FK151,ZY23,ZY24},  nullnorms and so on are introduced  on the unit interval $[0, 1]$ and then  can be  generalized  to bounded lattices.
 Considering the construction of these aggregation functions,  many researchers study them on the  bounded lattices.  Compared with  $[0,1]$,  bounded lattices are not totally-ordered. That is, the total order is a non-essential condition for the construction methods for aggregation functions. Similarly,  Compared with  bounded lattices,  the order on a  trellis does not need to satisfy transitivity.
  In this case, the concepts of $t$-norms and uninorms had been generalized from the bounded lattices to bounded trellises by  Zedam, De Baets \cite{LZ23} and Kong, Zhao \cite{YK24}, respectively. This shows that the  transitivity of all elements  is a non-essential condition for the construction methods for these aggregation functions on bounded trellises. So, the concepts of transitive elements, left-, right- and   middle-transitive  elements are introduced  \cite{LZ23}. In other words,  the construction methods for these aggregation functions do just need the condition that some of the  elements
 are transitive, left-, right- or middle-transitive.

As in \cite{LZ23},  Zedam and  De Baets wrote: our goal is to unravel how abandoning the transitivity property affects the notion  of a $t$-norm on bounded psosets, and in particular on bounded trellises. So, we also want to know  how to  construct   nullnorms on bounded trellises. Specifically, we discuss what the role  the transitivity  play in the construction  for nullnorms on bounded trellises.
In this paper, we introduce the notion of nullnorms on bounded trellises and propose some construction methods for nullnorms on bounded trellises. Moreover, we find that the zero element of nullnorms on bounded trellises must be middle-transitive.

%
%
%
%
%
%

\section{Preliminaries}
In this section, we recall some basic notions and results related to trellises.

\begin{definition}[\cite{EF70}]\label{}
Let $X$ be a set. A pseudo order on $X$ is a binary relation $\unlhd$ on $X$ satisfies:\\
(i) reflexive: $x\unlhd x$, for any $x\in X$;\\
(ii) antisymmetric: $x\unlhd y$ and $y\unlhd x$ implies $x=y$, for any $x,y\in X$.

A set $X$ equipped with a pseudo-order $\unlhd$ is called a pseudo-ordered set (psoset, for short) and is denoted by $\mathbb{P}=(X,\unlhd)$.
\end{definition}

Let $\mathbb{P}=(X,\unlhd)$ be a psoset, $a,b\in X$  and $A\subseteq X$. The symbol $a\lhd b$ means that $a\unlhd b$ and $a\neq b$. If neither $a\unlhd b$ nor $b\unlhd a$, then we say that $a$ and $b$ are incomparable and we use the notation
$a\parallel b$. The set of all elements of $X$ that are incomparable with $a$ is denoted by $I_{a}$, i.e., $I_{a}=\{x\in X\mid x\parallel a\}$.
 The antisymmetry of the pseudo-order implies that if $A$ has an infimum (resp. supremum), then it is unique, and is denoted by $\bigwedge A$ (resp. $\bigvee A$). Specially, if $A=\{a,b\}$, then we write $a\wedge b$ (called meet) instead of $\bigwedge\{a,b\}$ and $a\vee b$ (called join) instead of $\bigvee\{a,b\}$ \cite{HS71}.


Let $\mathbb{P}=(X,\unlhd)$ be a psoset, $\mathcal{C}$ be a subset of $X$, $x,y\in X$ and $z,t\in \mathcal{C}$. The symbol $x\lesssim y$ means that there exists a finite sequence
$(x_{1},\ldots,x_{n})$ of elements from $X$ such that $x\unlhd x_{1}\unlhd \ldots\unlhd x_{n}\unlhd y$. The symbol $z\lesssim_{\mathcal{C}} t$ means that there exists a finite sequence $(c_{1},\dots, c_{n})$ of elements from $\mathcal{C}$ such that $z\unlhd c_{1}\unlhd \ldots\unlhd c_{n}\unlhd t$. If for any $z,t\in\mathcal{C}$, both $z\lesssim_{\mathcal{C}} t$ and $t\lesssim_{\mathcal{C}} z$ hold, then $\mathcal{C}$ is called a cycle.  Due to the antisymmetry of $\unlhd$, any non-trivial cycle contains at least three elements.

Similarly,  as for partially ordered sets (posets, for short), a finite pseudo-ordered set can be represented by a Hasse-type diagram (or, simply, Hasse diagram) with the following difference: if $x$ and $y$ are not related, while in a poset set this would be implied by transitivity, then $x$ and $y$ are joined by a dashed edge. If $x\lesssim y$ and $y\unlhd x$, then $x$ and $y$ are joined by a directed edge going from $y$ to $x$ \cite{LZ23}.

\begin{example}\label{ex21}
Let $\mathbb{P}=(\{a, b, c, d\}, \unlhd)$ be a psoset given by Hasse diagram in Fig.1. Here,
$b\unlhd c$, $c\unlhd d$, while $b\ntrianglelefteq d$. And $\{a, b, c, d\}$ is a cycle.
\end{example}

\begin{figure}
 \centering
 \includegraphics[width=4cm, height=4cm]{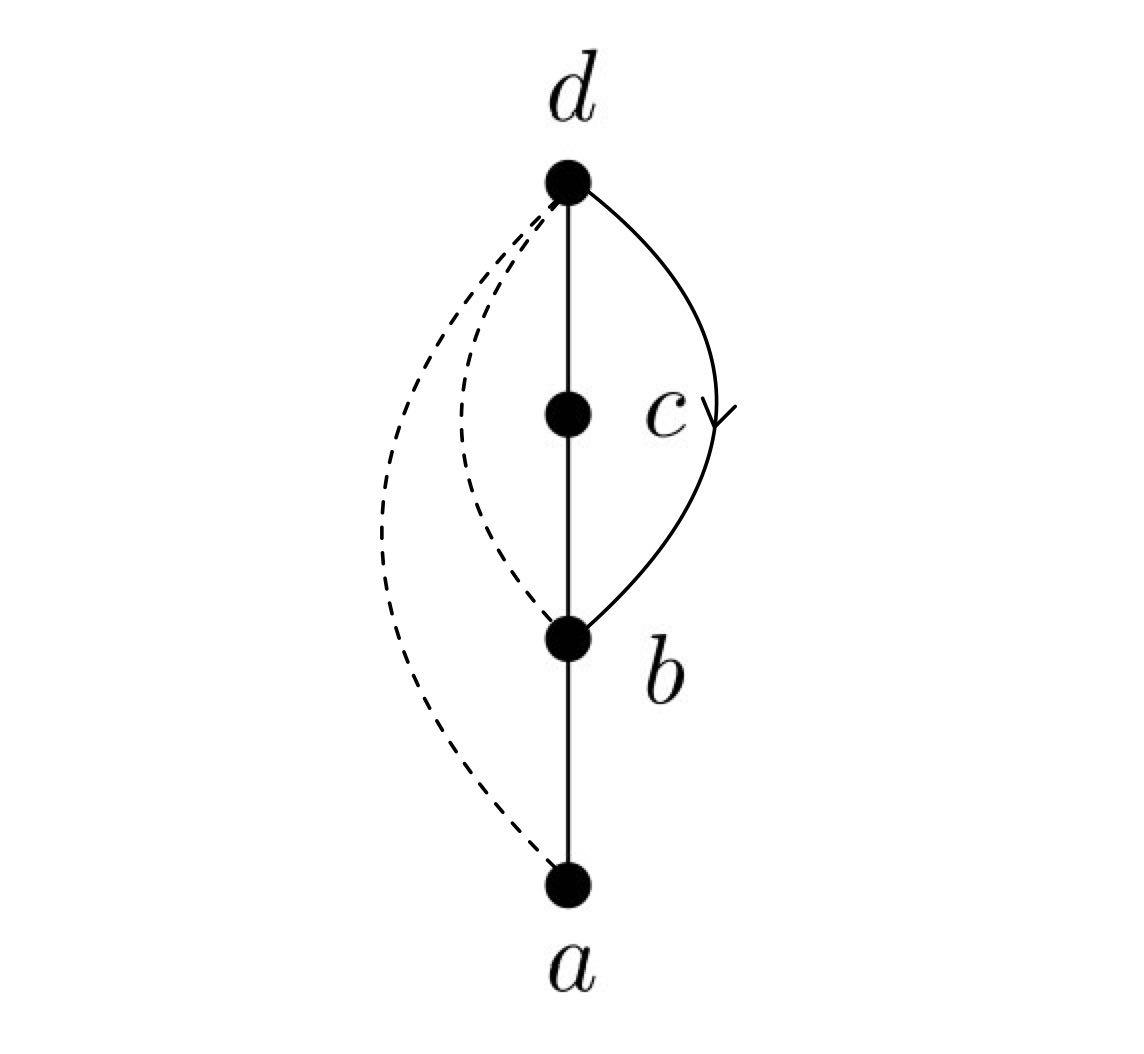}

 $\put(-130,0){\emph{Fig.1. Hasse diagram of the psoset in Example \ref{ex21}.}}$
\end{figure}

%
%
%
%
%
%
%
%
%
%
%
%
%

\begin{definition}[\cite{YK24}]\label{}
A psoset $\mathbb{P}=(X,\unlhd)$ is called bounded if it has a
smallest element denoted by $0$ and a greatest element denoted by $1$, i.e., $0\unlhd x\unlhd 1$,
for any $x\in X$, and is denoted by $\mathbb{P}=(X,\unlhd,0,1)$.
\end{definition}

\begin{definition}[\cite{KG73}]\label{}
A $\wedge$-semi-trellis (resp. $\vee$-semi-trellis) is a psoset $\mathbb{P}=(X,\unlhd)$ such that $x\wedge y$
(resp. $x\vee y$) exists, for any $x,y\in X$. A trellis is a psoset that is both a $\wedge$-semi-trellis and a $\vee$-semi-trellis; it is denoted by $\mathbb{P}=(X,\unlhd,\wedge,\vee)$. It is clear that $a\unlhd b$ is defined as $a\wedge b=a$/$a\vee b=b$ for a
trellis. A bounded trellis is denoted by $\mathbb{P}=(X,\unlhd,\wedge,\vee,0,1)$. In this paper, we consider
the bounded trellis with at least three elements.
\end{definition}

\begin{definition}[\cite{YK24}]\label{}
Let $\mathbb{P}=(X,\unlhd,\wedge,\vee,0,1)$ be a bounded trellises, $a,b\in X$ with $a\unlhd b$. A subinterval $[a,b]$ of $X$ is defined as
$$[a,b]=\{x\in X\mid a\unlhd x \unlhd b\}.$$
Similarly, we can define $[a,b[=\{x\in X\mid a\unlhd x \lhd b\}, ]a,b]=\{x\in L\mid a\lhd x \unlhd b\}$ and $]a,b[=\{x\in L\mid a\lhd x \lhd b\}$.
\end{definition}

For any $a\in X$, let $I_{a}^{1}=\{x\in I_{a}\mid x\lesssim a$ and $x$ is joined to $a$ by a dashed edge
$\}$, $I_{a}^{2}=\{x\in I_{a}\mid a\lesssim x$ and $a$ is joined to $x$ by a dashed edge $\}$, $I_{a}^{3}=\{x\in I_{a}\mid x$ and $a$ are not joined by a dashed edge $\}$,
$N(a)=\{x\in I_{a}^{3}\mid \exists y\in[0,a[\cup I_{a}^{1}$ such that $x\unlhd y\}$, $M(a)=\{x\in I_{a}^{3}\mid \exists y\in]a,1]\cup I_{a}^{2}$ such that $y\unlhd x\}$, $N_{i}=\{x\in\ ]0,a[\ \mid \exists\ y\in I_{a}^{1}$ such that $y\lhd x\}$, $M_{i}=\{x\in\ ]a,1[\ \mid \exists\ y\in I_{a}^{2}$ such that $x\lhd y\}$ and $K=\{x\in X^{mtr}\mid x$ does not belong to any cycle$\}$.

\begin{definition}[\cite{HS72}]\label{}
Let $\mathbb{P}=(X,\unlhd)$ be a psoset. An element $a\in X$ is called:\\
(i) right-transitive, if $a\unlhd x\unlhd y$ implies $a\unlhd y$, for any $x,y\in X$. We denote by $X^{rtr}$ the set of right-transitive elements of $X$;\\
(ii) left-transitive, if $x\unlhd y\unlhd a$ implies $x\unlhd a$, for any $x,y\in X$. We denote by $X^{ltr}$ the set of left-transitive elements of $X$;\\
(iii) middle-transitive, if $x\unlhd a\unlhd y$ implies$x\unlhd y$, for any $x,y\in X$. We denote by $X^{mtr}$ the set of middle-transitive elements of $X$;\\
(iv) transitive, if it is right-, left- and middle-transitive. We denote by $X^{tr}$ the set of transitive elements of $X$.
\end{definition}

\begin{remark}[\cite{YK24}]
If $a\in X^{mtr}$, then $x\unlhd y$ for any $x\in[0,a]$ and $y\in[a,1]$ due to $x\unlhd a\unlhd y$.
\end{remark}

\begin{proposition}[\cite{HS72}]\label{}
Let $\mathbb{P}=(X,\unlhd,\wedge,\vee)$ be a trellis and $a\in X$.\\
(i) If $a$ is right-transitive, then $a\unlhd x$ implies $a\vee y\unlhd x\vee y$, for any $x,y\in X$.\\
(ii) If $a$ is left-transitive, then $x\unlhd a$ implies $x\wedge y\unlhd a\wedge y$, for any $x,y\in X$.
\end{proposition}

\begin{definition}[\cite{LZ23}]\label{}
Let $\mathbb{P}=(X,\unlhd)$ be a psoset. A binary operation $F$ on $\mathbb{P}$ is called:\\
(i) commutative, if $F(x,y)=F(y,x)$, for any $x,y\in X$;\\
(ii) associative, if $F(x,F(y,z))=F(F(x,y),z)$, for any $x,y,z\in X$;\\
(iii) right-increasing, if $x\unlhd y$ implies $F(z,x)\unlhd F(z,y)$, for any $x,y,z\in X$;\\
(iv) left-increasing, if $x\unlhd y$ implies $F(x,z)\unlhd F(y,z)$, for any $x,y,z\in X$;\\
(v) increasing, if $x\unlhd y$ and $z\unlhd t$ implies $F(x,z)\unlhd F(y,t)$, for any $x,y,z,t\in X$.
\end{definition}

\begin{definition}[\cite{LZ23}]\label{}
Let $\mathbb{P}=(X,\unlhd,0,1)$ be a bounded psoset. A binary operation $T: X^{2}\rightarrow X$ ($S: X^{2}\rightarrow X$) is called a triangular norm (triangular conorm) on $\mathbb{P}$ if it is increasing, commutative,
associative and has $1$ as neutral element ($0$ as neutral element), i.e., $T(x,1)=x$ ($S(x,0)=x$), for any
$x\in X$.
\end{definition}

\section{Nullnorms on bounded trellises }
In this section, we introduce the concept of nullnorms on bounded trellises. And we give some remarks with respect to the zero
element and the structure of nullnorms. Moreover, we provide some methods of constructing nullnorms on bounded trellises.

\begin{definition}
Let $\mathbb{P}=(X,\unlhd,\wedge,\vee,0,1)$ be a bounded trellis. An operation $V:X^{2}\rightarrow X$ is called a nullnorm on $\mathbb{P}$ if it is commutative, associative, increasing and has an element $a\in X$ such that $V(x,0)=x$ for all $x\unlhd a$, $V(x,1)=x$ for all $x\unrhd a$. If $a\in X\setminus\{0, 1\}$, then $V$ is called a proper nullnorm.
\end{definition}

%
%
%

\begin{proposition}\label{pro31}
Let $\mathbb{P}=(X,\unlhd,\wedge,\vee,0,1)$ be a bounded trellis, $V$ be a nullnorm on $\mathbb{P}$ with zero
element $a$ and $x\in X$. \\
$(1)$ If $x\unlhd a$, then $V(x,y)=a$ for any $y\in [a,1]$.\\
$(2)$ If $a\unlhd x$, then $V(x,y)=a$ for any $y\in [0,a]$.
\end{proposition}
\begin{proof}
$(1)$ Let $x\in X$ with $x\unlhd a$, and $y\in [a,1]$. Then $a=V(0,a)\unlhd V(x,y)$ and $V(x,y)\unlhd V(a,1)=a$. Thus, $V(x,y)=a$.\\
$(2)$ Let $x\in X$ with $a\unlhd x$, and $y\in [0,a]$. Then $a=V(a,0)\unlhd V(x,y)$ and $V(x,y)\unlhd V(1,a)=a$. Thus, $V(x,y)=a$.
\end{proof}

\begin{proposition}
Let $\mathbb{P}=(X,\unlhd,\wedge,\vee,0,1)$ be a bounded trellis and $V$ be a nullnorm on $\mathbb{P}$ with zero
element $a$ and $x\in X$. \\
$(1)$ If $x\in I_{a}^{1}$, then $V(x,y)=a$ for any $y\in [a,1]$.\\
$(2)$ If $x\in I_{a}^{2}$, then $V(x,y)=a$ for any $y\in [0,a]$.
\end{proposition}
\begin{proof}
(1) Let $x\in I_{a}^{1}$. Then there exist a finite sequence $(x_{1},\dots,x_{n},y_{1},\dots,y_{n},z_{1},\dots,z_{n})$  such that $x\unlhd x_{1}\unlhd \dots\unlhd x_{n}\unlhd y_{1}\unlhd \dots\unlhd y_{n}\unlhd z_{1}\unlhd \dots\unlhd z_{n}\unlhd a$,  where $(x_{1},\dots,x_{n})\subseteq I_{a}^{1}$,  $(y_{1},\dots,y_{n})\subseteq X$  and  $(z_{1},\dots,z_{n})\subseteq [0,a]$. Take $y\in [a,1]$. Then we can obtain the following results:

Since $z_{n}\unlhd a$, $V(z_{n},y)=a$ by Proposition \ref{pro31}(1).

Since $z_{n-1}\unlhd z_{n}$, $a=V(0,a)\unlhd V(z_{n-1},y)\unlhd V(z_{n},y)=a$. Thus, $V(z_{n-1},y)=a$.

Since $z_{n-2}\unlhd z_{n-1}$, $a=V(0,a)\unlhd V(z_{n-2},y)\unlhd V(z_{n-1},y)=a$. Thus, $V(z_{n-2},y)=a$.

By induction, we can obtain $V(x_{1},y)=a$.

Since $x\unlhd x_{1}$, $a=V(0,a)\unlhd V(x,y)\unlhd V(x_{1},y)=a$. Thus, $V(x,y)=a$.

Therefore, if $x\in I_{a}^{1}$, then $V(x,y)=a$ for any $y\in [a,1]$.

(2) This can be proved similarly as (1).
\end{proof}

\begin{proposition}
Let $\mathbb{P}=(X,\unlhd,\wedge,\vee,0,1)$ be a bounded trellis and $V$ be a nullnorm on $\mathbb{P}$ with zero
element $a$. Then:\\
$(1)$ $V\mid_{[0,a]^{2}}:[0,a]^{2}\rightarrow[0,a]$ is a $t$-conorm on $[0,a]$.\\
$(2)$ $V\mid_{[a,1]^{2}}:[a,1]^{2}\rightarrow[a,1]$ is a $t$-norm on $[a,1]$.\\
$(3)$ $V(x,y)=a$, for all $(x,y)\in[0,a]\times[a,1]\cup [a,1]\times[0,a]$.\\
$(4)$ $a\unlhd V(x,y)$, for all $(x,y)\in[a,1]^{2}\cup [a,1]\times I_{a}\cup I_{a}\times[a,1]$.\\
$(5)$ $V(x,y)\unlhd a$, for all $(x,y)\in[0,a]^{2}\cup [0,a]\times I_{a}\cup I_{a}\times[0,a]$.\\
$(6)$ $V(x,y)\unlhd y$, for all $(x,y)\in X\times[a,1]$.\\
$(7)$ $V(x,y)\unlhd x$, for all $(x,y)\in[a,1]\times X$.\\
$(8)$ $x\unlhd V(x,y)$, for all $(x,y)\in[0,a]\times X$.\\
$(9)$ $y\unlhd V(x,y)$, for all $(x,y)\in X\times[0,a]$.\\
$(10)$ $x\vee y\unlhd V(x,y)$, for all $(x,y)\in [0,a]^{2}$.\\
$(11)$ $V(x,y)\unlhd x\wedge y$, for all $(x,y)\in[a,1]^{2}$.\\
$(12)$ $(x\wedge a)\vee(y\wedge a)\unlhd V(x,y)$, for all $(x,y)\in [0,a]\times I_{a}\cup I_{a}\times[0,a]\cup I_{a}\times I_{a}$.\\
$(13)$ $V(x,y)\unlhd (x\vee a)\wedge(y\vee a)$, for all $(x,y)\in [a,1]\times I_{a}\cup I_{a}\times[a,1]\cup I_{a}\times I_{a}$.
\end{proposition}
\begin{proof}
Obviously, (1) and (2) hold.\\
(3) Since $V(x,y)\unlhd V(a,1)=a$ and $a=V(0,a)\unlhd V(x,y)$, $V(x,y)=a$. The other case can be proved immediately by the commutativity property of $V$.\\
(4) Based on the increasingness of $V$, we can obtain that $a=V(0,a)\unlhd V(x,y)$. The other case can be proved immediately by the commutativity property of $V$.\\
(5) Based on the increasingness of $V$, we can obtain that $V(x,y)\unlhd V(a,1)=a$. The other case can be proved immediately by the commutativity property of $V$.\\
(6) 
$V(x,y)\unlhd V(1,y)=y$, for all $(x,y)\in X\times[a,1]$.\\
(7) 
$V(x,y)\unlhd V(x,1)=x$, for all $(x,y)\in [a,1]\times X$.\\
(8) 
$x=V(x,0)\unlhd V(x,y)$, for all $(x,y)\in[0,a]\times X$.\\
(9) 
$y=V(0,y)\unlhd V(x,y)$, for all $(x,y)\in X\times [0,a] $.\\
(10) Let $(x,y)\in [0,a]^{2}$. Since $x=V(x,0)\unlhd V(x,y)$ and  $y=V(0,y)\unlhd V(x,y)$, $x\vee y\unlhd V(x,y)$. \\
(11) Let $(x,y)\in [a,1]^{2}$. Since $V(x,y)\unlhd V(x,1)=x$ and  $V(x,y)\unlhd V(1,y)=y$, $V(x,y)\unlhd x\wedge y$. \\
(12) Let $(x,y)\in [0,a]\times I_{a}\cup I_{a}\times [0,a]\cup  I_{a}\times I_{a}$. Since $x\wedge a=V(x\wedge a,0)\unlhd V(x,y)$ and $y\wedge a=V(0,y\wedge a)\unlhd V(x,y)$, $(x\wedge a)\vee(y\wedge a)\unlhd V(x,y)$.\\
(13) Let $(x,y)\in [a,1]\times I_{a}\cup I_{a}\times [a,1]\cup  I_{a}\times I_{a}$. Since $V(x,y)\unlhd V(x\vee a,1)=x\vee a$ and $V(x,y)\unlhd V(1,y\vee a)=y\vee a$, $V(x,y)\unlhd (x\vee a)\wedge(y\vee a)$.
\end{proof}

\begin{proposition}
Let $\mathbb{P}=(X,\unlhd,\wedge,\vee,0,1)$ be a bounded trellis and $V$ be a nullnorm on $\mathbb{P}$ with zero
element $a$. Then $a\in X^{mtr}$.
\end{proposition}
\begin{proof}
Let $a$ be the zero element of a nullnorm $V$ on $\mathbb{P}$. It is obvious that if at least one of $]0,a[$ and $]a,1[$ is an empty set, then $a\in X^{mtr}$. Next, we just give the proof with $]0,a[\neq \emptyset$ and $]a,1[\neq \emptyset$. Assume that $a\notin X^{mtr}$. Then there exist $b,c\in X$ such that $b\unlhd a$, $a\unlhd c$ and $b\ntrianglelefteq c$. Since $V$ is increasing, we have $b=V(b,0)\unlhd V(c,1)=c$.
However, this is a contradiction. Therefore, $a\in X^{mtr}$.
%
\end{proof}

Let $\mathbb{P}=(X,\unlhd,\wedge,\vee,0,1)$ be a bounded trellis.   We know that there exists a nullnorm $V$ for $a=0$ (In this case, $V$ is a $t$-norm) and
there exists a nullnorm $V$ for $a = 1$ (In this case, $V$ is a $t$-conorm). One can wonder whether there exist a nullnorm on
bounded trellises $\mathbb{P}$ for $a\in X\setminus\{0,1\}$. To illustrate the presence of such a nullnorm, we shall give the following theorem. By the following theorem, based on the presence of $t$-conorms on $[0, a]$ and the presence of $t$-norms on $[a, 1]$,  we can guarantee the existence of a nullnorm on $\mathbb{P}$ for some $a\in X\setminus\{0,1\}$.

\begin{theorem}\label{th31}
Let $\mathbb{P}=(X,\unlhd,\wedge,\vee,0,1)$ be a bounded trellis, $a\in K\setminus\{0,1\}$, $S:[0,a]^{2}\rightarrow[0,a]$ be a $t$-conorm and  $T:[a,1]^{2}\rightarrow[a,1]$ be a $t$-norm. Assume that $N_{i}\cup I_{a}^{1}\subseteq X^{ltr}\setminus\{0,1\}$ and $M_{i}\cup I_{a}^{2}\cup I_{a}^{3}\subseteq X^{rtr}\setminus\{0,1\}$. Then the binary operation  $V$ is a nullnorm on $\mathbb{P}$ with the zero element $a$ if and only if  $N(a)=\emptyset$, where $V:X^{2}\rightarrow X$ is defined by

$V(x,y)=\begin{cases}
S(x\wedge a, y\wedge a) &\mbox{if } (x,y)\in ([0,a[\cup I_{a}^{1})^{2},\\
T(x\vee a, y\vee a) &\mbox{if } (x,y)\in (]a,1]\cup I_{a}^{2}\cup I_{a}^{3})^{2},\\
a &\mbox{}otherwise.
\end{cases}$
\end{theorem}
\begin{proof}
Necessity:  Let $V$ be a nullnorm on $\mathbb{P}$ with the zero element $a$. We prove that $N_{a}=\emptyset$.

Assume that $N_{a}\neq\emptyset$. Take $x\in N_{a}$. Then there exists an element $y\in[0,a[\cup I_{a}^{1}$ such that $x\unlhd y$. We can obtain that $V(x,0)=a$ and $V(y,0)=S(y\wedge a,0\wedge a)=S(y\wedge a,0)=y\wedge a$. Since $a\ntrianglelefteq y\wedge a$, this contradicts the increasingness property of $V$. Thus $N_{a}=\emptyset$.

Sufficiency: We have that $V(x,0)=S(x,0)=x$ for all $x\leq a$ and $V(x,1)=T(x,1)=x$ for all $x\geq a$. So, $a$ is the zero element of $V$. It is obvious that $V$ is commutative. Hence, we give only the proof of the increasingness and the associativity of $V$.

Let $x,y,z,t\in X$ with $x\unlhd y$, $z\unlhd t$. We divided the following cases to prove $V(x,z)\unlhd V(y,t)$:

(1) $x,y,z,t=a$. Then $V(x,z)=a=V(y,t)$.

(2) $x,y,z=a$, $t\neq a$. Then $V(x,z)=a=V(y,t)$. The proof of $x,y,t=a$, $z\neq a$; $x,z,t=a$, $y\neq a$; $y,z,t=a$, $x\neq a$ is similar.

(3) $x,t=a$, $z\neq a$, $y\neq a$. Then $V(x,z)=a=V(y,t)$.
The proof of $y,z=a$, $x\neq a$, $t\neq a$; $x,y=a$, $z\neq a$, $t\neq a$; $z,t=a$, $x\neq a$, $y\neq a$ is similar.

(4) $y,t=a$, $z\lhd a$, $x\lhd a$. Then $V(x,z)=S(x\wedge a,z\wedge a)\unlhd S(a,a)=a=V(y,t)$.

(5) $x,z=a$, $a\lhd y$, $a\lhd t$. Then $V(x,z)=a=T(a,a)\unlhd T(y\vee a,t\vee a)=a=V(y,t)$.

(6) $y,z,t\neq a$, $x=a$. Then $V(x,z)=a$ and $a\lhd y$. If $t\in [0,a[\cup I_{a}^{1}$, then $V(x,z)=a=V(y,t)$.
If $t\in I_{a}^{2}\cup I_{a}^{3}\cup ]a,1]$, then $V(x,z)=a=T(a,a)\unlhd T(y\vee a,t\vee a)=V(y,t)$.

(7) $x,y,t\neq a$, $z=a$. Then $V(x,z)=a$ and $a\lhd t$. If $y\in [0,a[\cup I_{a}^{1}$, then $V(x,z)=a=V(y,t)$.
If $y\in I_{a}^{2}\cup I_{a}^{3}\cup ]a,1]$, then $V(x,z)=a=T(a,a)\unlhd T(y\vee a,t\vee a)=V(y,t)$.

(8) $x,z,t\neq a$, $y=a$. Then $V(y,t)=a$ and $x\lhd a$. If $z\in [0,a[\cup I_{a}^{1}$, then $V(x,z)=S(x\wedge a,z\wedge a)\unlhd S(a,a)=a=V(y,t)$.
If $z\in I_{a}^{2}\cup I_{a}^{3}\cup ]a,1]$, then $V(x,z)=a=V(y,t)$.

(9) $x,y,z\neq a$, $t=a$. Then $V(y,t)=a$ and $z\lhd a$. If $x\in [0,a[\cup I_{a}^{1}$, then $V(x,z)=S(x\wedge a,z\wedge a)\unlhd S(a,a)=a=V(y,t)$.
If $x\in I_{a}^{2}\cup I_{a}^{3}\cup ]a,1]$, then $V(x,z)=a=V(y,t)$.

(10) $x,y,z,t\neq a$.

(10-a) If $x,z\in [0,a[\cup I_{a}^{1}$, then $V(x,z)=S(x\wedge a,z\wedge a)$.

Whenever $y,t\in [0,a[\cup I_{a}^{1}$, then $V(y,t)=S(y\wedge a,t\wedge a)$. Since $M_{i}\cup I_{a}^{1}\subseteq X^{ltr}\setminus\{0,1\}$, we have that $x\wedge a\unlhd y\wedge a$ and $z\wedge a\unlhd t\wedge a$, and then $V(x,z)\unlhd V(y,t)$.

Whenever $y\notin [0,a[\cup I_{a}^{1}$, $t\in [0,a[\cup I_{a}^{1}$, then we have that $V(y,t)=a$ and then $V(x,z)=S(y\wedge a,t\wedge a)\unlhd S(a,a)=a=V(y,t)$.

Whenever $y\in [0,a[\cup I_{a}^{1}$, $t\notin [0,a[\cup I_{a}^{1}$, then we have that $V(y,t)=a$ and then $V(x,z)=S(y\wedge a,t\wedge a)\unlhd S(a,a)=a=V(y,t)$.

Whenever $y,t\notin [0,a[\cup I_{a}^{1}$, then we have that $V(y,t)=T(y\vee a,t\vee a)$ and then $V(x,z)=S(y\wedge a,t\wedge a)\unlhd S(a,a)=a=T(a,a)\unlhd T(y\vee a,t\vee a)=V(y,t)$.

(10-b) If $x\in [0,a[\cup I_{a}^{1}$,  $z\notin [0,a[\cup I_{a}^{1}$, then $t\notin [0,a[\cup I_{a}^{1}$ and $V(x,z)=a$.

Whenever $y\in [0,a[\cup I_{a}^{1}$, then we have that $V(y,t)=a$ and then $V(x,z)=a=V(y,t)$.

Whenever $y\notin [0,a[\cup I_{a}^{1}$, then we have that $V(y,t)=T(y\vee a,t\vee a)$ and then $V(x,z)=a=T(a,a)\unlhd T(y\vee a,t\vee a)=V(y,t)$.

(10-c) If $x\notin [0,a[\cup I_{a}^{1}$,  $z\in [0,a[\cup I_{a}^{1}$, then $y\notin [0,a[\cup I_{a}^{1}$ and $V(x,z)=a$.

Whenever $t\in [0,a[\cup I_{a}^{1}$, then we have that $V(y,t)=a$ and then $V(x,z)=a=V(y,t)$.

Whenever $t\notin [0,a[\cup I_{a}^{1}$, then we have that $V(y,t)=T(y\vee a,t\vee a)$ and then $V(x,z)=a=T(a,a)\unlhd T(y\vee a,t\vee a)=V(y,t)$.

(10-d) If $x,z\notin [0,a[\cup I_{a}^{1}$, then $y,t\notin [0,a[\cup I_{a}^{1}$. Since $M_{ii}\cup I_{a}^{2}\cup I_{a}^{3}\subseteq X^{rtr}\setminus\{0,1\}$, we have that $x\wedge a\unlhd y\wedge a$ and $z\wedge a\unlhd t\wedge a$, and then $V(x,z)=T(x\vee a,z\vee a)\unlhd T(y\vee a,t\vee a)=V(y,t)$.

Therefore, $V$ is increasing.  Let $x,y,z\in X$, we divided the four cases to prove $V(V(x,y),z)=V(x,V(y,z))$:

(1) $x,y,z=a$. Then $V(V(x,y),z)=a=V(x,V(y,z))$.

(2) $x,y=a$, $z\neq a$. Then $V(V(x,y),z)=a=V(x,V(y,z))$. The proof of $x,z=a$, $y\neq a$ and $y,z=a$, $x\neq a$ is similar.

(3) $x,y\neq a$, $z=a$. Then $V(V(x,y),z)=a=V(x,V(y,z))$. The proof of $x,z\neq a$, $y=a$ and $y,z\neq a$, $x=a$ is similar.

(4) $x,y,z\neq a$.

(4-a) If $x,y\in [0,a[\cup I_{a}^{1}$ and $z\notin [0,a[\cup I_{a}^{1}$, then $V(V(x,y),z)=V(S(x\wedge a,y\wedge a),z)=a=V(x,V(y,z))$.

(4-b) If $x,z\in [0,a[\cup I_{a}^{1}$ and $y\notin [0,a[\cup I_{a}^{1}$, then $V(V(x,y),z)=V(a,z)=a=V(x,a)=V(x,V(y,z))$.

(4-c) If $y,z\in [0,a[\cup I_{a}^{1}$ and $x\notin [0,a[\cup I_{a}^{1}$, then $V(V(x,y),z)=V(a,z)=a=V(x,S(y\wedge a,z\wedge a))=V(x,V(y,z))$.

(4-d) If $x\in [0,a[\cup I_{a}^{1}$ and $y,z\notin [0,a[\cup I_{a}^{1}$, then $V(V(x,y),z)=V(a,z)=a=V(x,T(y\vee a,z\vee a))=V(x,V(y,z))$.

(4-e) If $y\in [0,a[\cup I_{a}^{1}$ and $x,z\notin [0,a[\cup I_{a}^{1}$, then $V(V(x,y),z)=V(a,z)=a=V(x,a)=V(x,V(y,z))$.

(4-f) If $z\in [0,a[\cup I_{a}^{1}$ and $x,y\notin [0,a[\cup I_{a}^{1}$, then $V(V(x,y),z)=V(T(x\vee a,y\vee a),z)=a=V(x,a)=V(x,V(y,z))$.

(4-g) If $x,y,z\notin [0,a[\cup I_{a}^{1}$, then $V(V(x,y),z)=V(T(x\vee a,y\vee a),z)=T(T(x\vee a,y\vee a)\vee a,z\vee a)=T(T(x\vee a,y\vee a),z\vee a)=T(x\vee a,T(y\vee a,t\vee a))=T(x\vee a,T(y\vee a,t\vee a)\vee a)=V(x,T(y\vee a,t\vee a))=V(x,V(y,z))$.

Therefore, $V$ is a nullnorm on $\mathbb{P}$ with the zero element $a$.
\end{proof}\

The Fig. 2 gives the structure of $V$ in Theorem \ref{th31}.

\begin{minipage}{11pc}
\setlength{\unitlength}{0.75pt}\begin{picture}(600,325)
\put(115,30){\makebox(0,0)[l]{\footnotesize$0$}}
\put(221,30){\makebox(0,0)[l]{\footnotesize$a$}}
\put(255,30){\makebox(0,0)[l]{\footnotesize$1$}}
\put(265,30){\makebox(0,0)[l]{\footnotesize$I_{a}^{2}\cup I_{a}^{3}$}}
\put(362,30){\makebox(0,0)[l]{\footnotesize$I_{a}^{1}$}}
\put(114,135){\makebox(0,0)[l]{\footnotesize$a$}}
\put(115,180){\makebox(0,0)[l]{\footnotesize$1$}}
\put(85,200){\makebox(0,0)[l]{\footnotesize$I_{a}^{2}\cup I_{a}^{3}$}}
\put(112,270){\makebox(0,0)[l]{\footnotesize$I_{a}^{1}$}}

\put(138,88){\makebox(0,0)[l]{\footnotesize$S(x\wedge a, y\wedge a)$}}
\put(177,178){\makebox(0,0)[l]{\footnotesize$a$}}
\put(138,268){\makebox(0,0)[l]{\footnotesize$S(x\wedge a, y\wedge a)$}}
\put(268,88){\makebox(0,0)[l]{\footnotesize$a$}}
\put(228,178){\makebox(0,0)[l]{\footnotesize$T(x\vee a, y\vee a)$}}
\put(268,268){\makebox(0,0)[l]{\footnotesize$a$}}
\put(318,88){\makebox(0,0)[l]{\footnotesize$S(x\wedge a, y\wedge a)$}}
\put(357,178){\makebox(0,0)[l]{\footnotesize$a$}}
\put(318,268){\makebox(0,0)[l]{\footnotesize$S(x\wedge a, y\wedge a)$}}

\put(270,45){\line(1,0){135}}
\put(270,45){\line(-1,0){135}}
\put(270,135){\line(1,0){135}}
\put(270,135){\line(-1,0){135}}
\put(270,225){\line(1,0){135}}
\put(270,225){\line(-1,0){135}}
\put(270,315){\line(1,0){135}}
\put(270,315){\line(-1,0){135}}
\put(225,45){\line(0,1){270}}
\put(315,45){\line(0,1){270}}
\put(135,45){\line(0,1){270}}
\put(405,45){\line(0,1){270}}

\put(135,0){\emph{Fig.2. The nullnorm $V$ in Theorem \ref{th31}.}}
\end{picture}\\
\end{minipage}

If $\mathbb{P}$ is a bounded lattice, then we can obtain the nullnorm $V_{S}^{T}$ on bounded lattices which is constructed by  Ertu\u{g}rul (\cite{UE18}, Theorem 6).

\begin{corollary}
Let $\mathbb{P}$ be a bounded lattice, $a\in \mathbb{P}\setminus\{0,1\}$, $S:[0,a]^{2}\rightarrow[0,a]$ be a $t$-conorm and $T:[a,1]^{2}\rightarrow[a,1]$ be a $t$-norm. Then the binary operation  $V_{S}^{T}$ is a nullnorm on $\mathbb{P}$ with the zero element $a$, where $V_{S}^{T}:X^{2}\rightarrow X$ is defined by

$V_{S}^{T}(x,y)=\begin{cases}
S(x,y) &\mbox{if } (x,y)\in [0,a[^{2},\\
T(x\vee a, y\vee a) &\mbox{if } (x,y)\in (]a,1]\cup I_{a})^{2},\\
a &\mbox{}otherwise.
\end{cases}$
\end{corollary}

The next example illustrates the construction method of nullnorms on bounded trellises in Theorem \ref{th31}.

\begin{figure}\label{fig3}
 \centering
 \includegraphics[width=5cm, height=8cm]{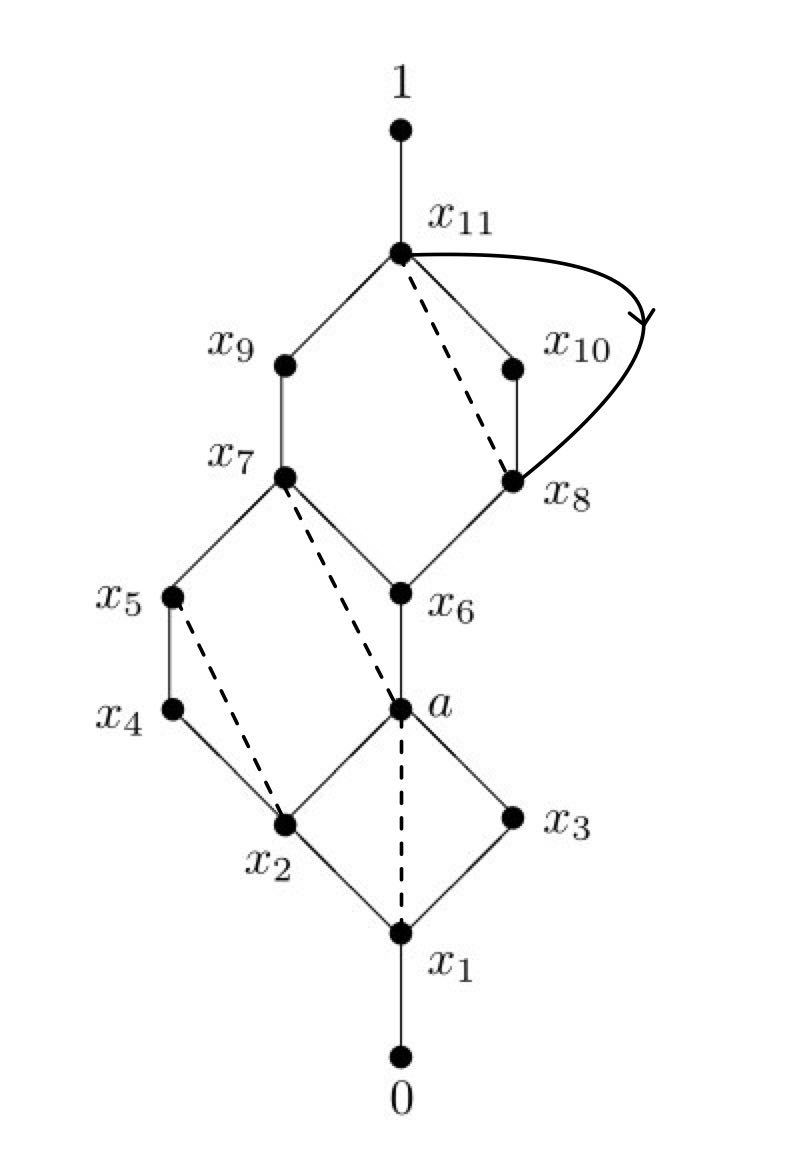}

 $\put(-130,0){\emph{Fig.3. Hasse diagram of the psoset in Example \ref{ex31}.}}$
\end{figure}

\begin{example}\label{ex31}
Let $\mathbb{P}=(\{0, a, x_{1}, x_{2}, x_{3}, x_{4}, x_{5}, x_{6}, x_{7}, x_{8}, x_{9}, 1\}, \unlhd, \wedge, \vee, 0, 1)$ be a bounded trellis
given by Hasse diagram in Fig. \ref{fig3}. Then $[0,a[=\{0,x_{2}, x_{3}\}$, $I_{a}^{1}=\{x_{1}\}$, $I_{a}^{2}=\{x_{7}\}$, $I_{a}^{3}=\{x_{4},x_{5}\}$, $]a,1]=\{x_{6}, x_{7}, x_{8}, x_{9}, x_{10}, x_{11}\}$, $N_{i}=\{x_{2}, x_{3}\}$ and $M_{i}=\{x_{6}\}$. It is obvious that $a\in K\setminus\{0,1\}$, $N_{i}\cup I_{a}^{1}\subseteq X^{ltr}\setminus\{0,1\}$, $M_{i}\cup I_{a}^{2}\cup I_{a}^{3}\subseteq X^{rtr}\setminus\{0,1\}$ and $N(a)=\emptyset$.
Given a $t$-conorm $S$ on $[0,a]$ defined by
$S(x,y)=\begin{cases}
x\vee y &\mbox{if x=0 or y=0,}\\
a &\mbox{}otherwise.
\end{cases}$
and a $t$-norm $T$ on $[a,1]$ defined by
$T(x,y)=\begin{cases}
x\wedge y &\mbox{if x=1 or y=1,}\\
a &\mbox{}otherwise.
\end{cases}$
(\cite{YK24}, Remark 3.4).
By using the construction method in Theorem \ref{th31}, we can obtain a nullnorm $V$ with the zero element $a$ shown in Table \ref{Tab:01}.
\end{example}

\begin{table}[htbp]
\centering
\caption{The nullnorm $V $ constructed by Theorem \ref{th31}.}
\label{Tab:01}\

\begin{tabular}{c|c c c c c c c c c c c c c c}
\hline
  $V$ & $0$ & $x_{1}$ & $x_{2}$ & $x_{3}$ & $a$ & $x_{4}$ & $x_{5}$ & $x_{6}$ & $x_{7}$ & $x_{8}$ & $x_{9}$ & $x_{10}$ & $x_{11}$ & $1$ \\
\hline
  $0$ & $0$ & $0$ & $x_{2}$ & $x_{3}$ & $a$ & $a$ & $a$ & $a$ & $a$ & $a$ & $a$ & $a$ & $a$ & $a$ \\

  $x_{1}$ & $0$ & $0$ & $x_{2}$ & $x_{3}$ & $a$ & $a$ & $a$ & $a$ & $a$ & $a$ & $a$ & $a$ & $a$ & $a$ \\

  $x_{2}$ & $x_{2}$ & $x_{2}$ & $a$ & $a$ & $a$ & $a$ & $a$ & $a$ & $a$ & $a$ & $a$ & $a$ & $a$ & $a$ \\

  $x_{3}$ & $x_{3}$ & $x_{3}$ & $a$ & $a$ & $a$ & $a$ & $a$ & $a$ & $a$ & $a$ & $a$ & $a$ & $a$ & $a$  \\

  $a$ & $a$ & $a$ & $a$ & $a$ & $a$ & $a$ & $a$ & $a$ & $a$ & $a$ & $a$ & $a$ & $a$ & $a$ \\

  $x_{4}$ & $a$ & $a$ & $a$ & $a$ & $a$ & $a$ & $a$ & $a$ & $a$ & $a$ & $a$ & $a$ & $a$ & $x_{9}$ \\

  $x_{5}$ & $a$ & $a$ & $a$ & $a$ & $a$ & $a$ & $a$ & $a$ & $a$ & $a$ & $a$ & $a$ & $a$ & $x_{9}$ \\

  $x_{6}$ & $a$ & $a$ & $a$ & $a$ & $a$ & $a$ & $a$ & $a$ & $a$ & $a$ & $a$ & $a$ & $a$ & $x_{6}$ \\

  $x_{7}$ & $a$ & $a$ & $a$ & $a$ & $a$ & $a$ & $a$ & $a$ & $a$ & $a$ & $a$ & $a$ & $a$ & $x_{9}$ \\

  $x_{8}$ & $a$ & $a$ & $a$ & $a$ & $a$ & $a$ & $a$ & $a$ & $a$ & $a$ & $a$ & $a$ & $a$ & $x_{8}$ \\

  $x_{9}$ & $a$ & $a$ & $a$ & $a$ & $a$ & $a$ & $a$ & $a$ & $a$ & $a$ & $a$ & $a$ & $a$ & $x_{9}$ \\

  $x_{10}$ & $a$ & $a$ & $a$ & $a$ & $a$ & $a$ & $a$ & $a$ & $a$ & $a$ & $a$ & $a$ & $a$ & $x_{10}$ \\

  $x_{11}$ & $a$ & $a$ & $a$ & $a$ & $a$ & $a$ & $a$ & $a$ & $a$ & $a$ & $a$ & $a$ & $a$ & $x_{11}$ \\

  $1$ & $a$ & $a$ & $a$ & $a$ & $a$ & $x_{9}$ & $x_{9}$ & $x_{6}$ & $x_{9}$ & $x_{8}$ & $x_{9}$ & $x_{10}$ & $x_{11}$ & $1$ \\
\hline
\end{tabular}
\end{table}\

\begin{theorem}\label{th32}
Let $\mathbb{P}=(X,\unlhd,\wedge,\vee,0,1)$ be a bounded trellis, $a\in K\setminus\{0,1\}$, $S:[0,a]^{2}\rightarrow[0,a]$ be a $t$-conrm and  $T:[a,1]^{2}\rightarrow[a,1]$ be a $t$-norm. Assume that $N_{i}\cup I_{a}^{1}\cup I_{a}^{3}\subseteq X^{ltr}\setminus\{0,1\}$ and $M_{i}\cup I_{a}^{2}\subseteq X^{rtr}\setminus\{0,1\}$. Then the binary operation  $V$ is a nullnorm on $\mathbb{P}$ with the zero element $a$ if and only if $M_{a}=\emptyset$, where $V:X^{2}\rightarrow X$ is defined by

$V(x,y)=\begin{cases}
S(x\wedge a, y\wedge a) &\mbox{if } (x,y)\in ([0,a[\cup I_{a}^{1}\cup I_{a}^{3})^{2},\\
T(x\vee a, y\vee a) &\mbox{if } (x,y)\in (]a,1]\cup I_{a}^{2})^{2},\\
a &\mbox{}otherwise.
\end{cases}$
\end{theorem}
\begin{proof}
This can be proved similarly as Theorem \ref{th31}.
\end{proof}\

The structure of $V$ in Theorem \ref{th32} is shown in Fig. 4.

\begin{minipage}{11pc}
\setlength{\unitlength}{0.75pt}\begin{picture}(600,325)
\put(115,30){\makebox(0,0)[l]{\footnotesize$0$}}
\put(221,30){\makebox(0,0)[l]{\footnotesize$a$}}
\put(265,30){\makebox(0,0)[l]{\footnotesize$1$}}
\put(280,30){\makebox(0,0)[l]{\footnotesize$I_{a}^{2}$}}
\put(340,30){\makebox(0,0)[l]{\footnotesize$I_{a}^{1}\cup I_{a}^{3}$}}
\put(114,135){\makebox(0,0)[l]{\footnotesize$a$}}
\put(115,180){\makebox(0,0)[l]{\footnotesize$1$}}
\put(112,200){\makebox(0,0)[l]{\footnotesize$I_{a}^{2}$}}
\put(85,270){\makebox(0,0)[l]{\footnotesize$I_{a}^{1}\cup I_{a}^{3}$}}

\put(138,88){\makebox(0,0)[l]{\footnotesize$S(x\wedge a, y\wedge a)$}}
\put(177,178){\makebox(0,0)[l]{\footnotesize$a$}}
\put(138,268){\makebox(0,0)[l]{\footnotesize$S(x\wedge a, y\wedge a)$}}
\put(268,88){\makebox(0,0)[l]{\footnotesize$a$}}
\put(228,178){\makebox(0,0)[l]{\footnotesize$T(x\vee a, y\vee a)$}}
\put(268,268){\makebox(0,0)[l]{\footnotesize$a$}}
\put(318,88){\makebox(0,0)[l]{\footnotesize$S(x\wedge a, y\wedge a)$}}
\put(357,178){\makebox(0,0)[l]{\footnotesize$a$}}
\put(318,268){\makebox(0,0)[l]{\footnotesize$S(x\wedge a, y\wedge a)$}}

\put(270,45){\line(1,0){135}}
\put(270,45){\line(-1,0){135}}
\put(270,135){\line(1,0){135}}
\put(270,135){\line(-1,0){135}}
\put(270,225){\line(1,0){135}}
\put(270,225){\line(-1,0){135}}
\put(270,315){\line(1,0){135}}
\put(270,315){\line(-1,0){135}}
\put(225,45){\line(0,1){270}}
\put(315,45){\line(0,1){270}}
\put(135,45){\line(0,1){270}}
\put(405,45){\line(0,1){270}}

\put(135,0){\emph{Fig.4. The nullnorm $V$ in Theorem \ref{th32}.}}
\end{picture}\\
\end{minipage}

If $\mathbb{P}$ is a bounded lattice, then we can obtain the nullnorm $V_{T}^{S}$ on bounded lattices which is constructed by  Ertu\u{g}rul (\cite{UE18}, Theorem 6).

\begin{corollary}
Let $\mathbb{P}$ be a bounded lattice, $a\in \mathbb{P}\setminus\{0,1\}$, $S:[0,a]^{2}\rightarrow[0,a]$ be a $t$-conorm and  $T:[a,1]^{2}\rightarrow[a,1]$ be a $t$-norm. Then the binary operation  $V_{T}^{S}$ is a nullnorm on $\mathbb{P}$ with the zero element $a$, where $V_{T}^{S}:X^{2}\rightarrow X$ is defined by

$V_{T}^{S}(x,y)=\begin{cases}
S(x\wedge a, y\wedge a) &\mbox{if } (x,y)\in ([0,a[\cup I_{a})^{2},\\
T(x,y) &\mbox{if } (x,y)\in ]a,1]^{2},\\
a &\mbox{}otherwise.
\end{cases}$
\end{corollary}

\section{Conclusion}

In this paper, we introduce the concept nullnorms on bounded trellises and investigate the construction methods for nullnorms on bounded trellises.
About the results in this paper, we have listed the following remarks.

(1) Based on the existence of $t$-norms and $t$-conorms on arbitrary bounded trellises, we propose some construction methods for nullnorms on bounded trellises under additional constraints on the transitivity of some elements.  Moreover, the condition that  $N(a)=\emptyset$ is sufficient and necessary.

(2) The zero element of nullnorms on bounded trellises is just middle-transitive.

(3) Not all elements of bounded trellises need to satisfy the transitivity property when constructing nullnorms on bounded trellises.


In the future, we will consider more methods to construct nullnorms on bounded trellises and  investigate other aggregative operators on bounded trellises, such as uni-nullnorms and null-uninorms.

%


\begin{thebibliography}{99}

\bibitem{CA06} C. Alsina, M.J. Frank, B. Schweizer, Associative Functions. Triangular Norms and Copulas, World Scientific, Hackensack, 2006.

\bibitem{TC01} T. Calvo, B. De Baets, J. Fodor, The functional equations of Frank and Alsina for uninorms and nullnorms, Fuzzy Sets Syst. 120(2001)385--394.

\bibitem{GC18} G.D. \c{C}ayl{\i}, F. Kara\c{c}al, Idempotent nullnorms on bounded lattices, Inf. Sci. 425(2018)154--163.
\bibitem{GD181} G.D. \c{C}ayl{\i}, On a new class of t-norms and t-conorms on bounded lattices, Fuzzy Sets Syst. 332(2018)129--143.
\bibitem{GC19} G.D. \c{C}ayl{\i}, Some methods to obtain $t$-norms and $t$-conorms on bounded lattices, Kybernetika 55(2)(2019)273--294.
\bibitem{GC20} G.D. \c{C}ayl{\i}, Nullnorms on bounded lattices derived from $t$-norms and $t$-conorms, Inf. Sci. 512(2020)1134--1154.
\bibitem{GC201} G.D. \c{C}ayl{\i}, Construction methods for idempotent nullnorms on bounded lattices, Appl. Math. Comput. 366(2020).

\bibitem{DD85} D. Dubois, H. Prade, A review of fuzzy set aggregation connectives, Inf. Sci. 36(1985)85--121.
\bibitem{DD00} D. Dubois, H. Prade, Fundamentals of Fuzzy Sets, Kluwer Academic Publishers, Boston, 2000.
\bibitem{JD08} J. Drewniak, P. Dryga\'{s}, E. Rak, Distributivity between uninorms and nullnorms, Fuzzy Sets Syst. 159(2008)1646--1657.
\bibitem{PD04} P. Dryga\'{s}, A characterization of idempotent nullnorms, Fuzzy Sets Syst. 145 (2004) 455--461.
\bibitem{PD15} P. Dryga\'{s}, Distributivity between semi $t$-operators and semi nullnorms, Fuzzy Sets Syst. 264(2015)100--109.

\bibitem{BD06} B. De Baets, H. De Meyer, B. De Schuymer, S. Jenei, Cyclic evaluation of transitivity of reciprocal relations, Soc. Choice Welf. 26(2006)217--238.
\bibitem{BD15} B. De Baets, K. De Loof, H. De Meyer, A frequentist view on cycle-transitivity of reciprocal relations, Fuzzy Sets Syst. 281(2015)198--218.
\bibitem{UE18} \"{U}. Ertu\u{g}rul, Constructions of nullnorms on bounded lattices and equivalence relation on nullnorms, Fuzzy Sets Syst. 334(2018)94--109.

\bibitem{PF70} P.C. Fishburn, Intransitive indifference in preference theory: a survey, Oper. Res. 18(1970)207--228.
\bibitem{EF70} E. Fried, Tournaments and non-associative lattices, Ann. Univ. Sci. Bp. Rolando E\"{o}tv\"{o}s Nomin., Sect. Math. 13(1970)151--164.

\bibitem{KG73} K. Gladstien, A characterization of complete trellises of finite length, Algebra Univers. 3(1973)341--344.
\bibitem{MG09} M. Grabisch, J.L. Marichal, R. Mesiar, E. Pap, Aggregation Functions, Cambridge University Press, Cambridge, 2009.


\bibitem{XH22} X. Hua, New constructions of nullnorms on bounded lattices, Fuzzy Sets Syst. 439(2022)126--141.

\bibitem{BK02} B. Kerr, M.A. Riley, M.W. Feldman, B.J.M. Bohannan, Local dispersal promotes biodiversity in a real-life game of rock-paper-scissors, Nature 418(2002)171--174.
\bibitem{FK15} F. Kara\c{c}al, M.A. \.{I}nce, R. Mesiar, Nullnorms on bounded lattice, Inf. Sci. 325(2015)227--236.
\bibitem{FK151} F. Kara\c{c}al, R. Mesiar, Uninorms on bounded lattices, Fuzzy Sets Syst. 261(2015)33--43.
\bibitem{YK24} Y. Kong, B. Zhao, Uninorms on bounded trellises, Fuzzy Sets Syst. 481(2024)108898.


\bibitem{MM99} M. Mas, G. Mayor, J. Torrens, $t$-operators, Int. J. Uncertain. Fuzziness Knowl.-Based Syst. 7(1999)31--50.
\bibitem{MM02} M. Mas, G. Mayor, J. Torrens, The distributivity condition for uninorms and t-operators, Fuzzy Sets Syst. 128(2002)209--225.
\bibitem{RM98} R. Mesiar, E. Pap, Different interpretations of triangular norms and related operations, Fuzzy Sets Syst. 96(1998)183--189.


\bibitem{FQ05} F. Qin, B. Zhao, The distributive equations for idempotent uninorms and nullnorms, Fuzzy Sets Syst. 155(2005)446--458.

\bibitem{TR07} T. Reichenbach, M. Mobilia, E. Frey, Mobility promotes and jeopardizes biodiversity in rock-paper-scissors game, Nature 448(2007)1046--1049.

\bibitem{HS71} H. Skala, Trellis theory, Algebra Univers. 1(1971)218--233.
\bibitem{HS72} H. Skala, Trellis Theory, American Mathematical Soc., vol. 121, 1972.
\bibitem{XS20} X.R. Sun, H.W. Liu, Representation of nullnorms on bounded lattices, Inf. Sci. 539(2020)269--276.
\bibitem{SS06} S. Saminger, On ordinal sums of triangular norms on bounded lattices, Fuzzy Sets Syst. 157(10)(2006)1403--1416.

\bibitem{XW20} X. Wu, S. Liang, G.D. \c{C}ayl{\i}, Characterizing idempotent nullnorms on a special class of bounded lattices, Fuzzy Sets Syst. 427(2022)161--173.

\bibitem{ZY23} Z.Y. Xiu, X. Zheng, New construction methods of uninorms on bounded lattices via uninorms, Fuzzy Sets Syst. 465(2023)108535.
\bibitem{ZY24} Z.Y. Xiu, X. Zheng, A new approach to construct uninorms via uninorms on bounded lattices, Kybernetika, 60(2)(2024)125--149.


\bibitem{LZ23} L. Zedam, B. De Baets, Triangular norms on bounded trellises, Fuzzy Sets Syst. 462(2023)108468.

\end{thebibliography}
\end{document}